\documentclass{elsarticle}

\usepackage{amssymb,setspace,psfrag}
\usepackage[applemac]{inputenc}
\usepackage{setspace}
\usepackage{pgfplots}

\usepackage{amsmath,url}
\usepackage[titletoc]{appendix}
\usepackage{amssymb,widetext}
\usepackage{lscape,amstext,mathrsfs}
\usepackage{graphicx}
\usepackage{epstopdf}
\usepackage{tikz,pmat}
\usepackage{afterpage}
\usepackage{rotating}
\usepackage[graphicx]{realboxes}
\usepackage{pgfgantt}
\usepackage[europeanresistors,cuteinductors,emptydiodes]{circuitikz}
\usepackage{graphics,array} 
\usepackage{epsfig} 
\usepackage{epstopdf}
\usepackage{epstopdf}
\usepackage{mhchem,bm}
\usepackage{caption}
\usepackage{subcaption}
\usetikzlibrary{shapes,arrows}
\usetikzlibrary{arrows, decorations.markings}

\newcommand{\real}{{\mathbb R}}

\newtheorem{thm}{Theorem}[section]
\newtheorem{lem}[thm]{Lemma}
\newtheorem{prop}[thm]{Proposition}
\newtheorem{cor}[thm]{Corollary}

\newtheorem{defi}[thm]{Definition}
\newtheorem{remark}[thm]{Remark}
\newtheorem{problem}[thm]{Problem}

\newtheorem{exam}[thm]{Example}

\newenvironment{pf}{\paragraph{Proof:}}{\hfill$\square$}

\journal{Systems and Control Letters}

\begin{document}

\begin{frontmatter}

\title{Safety Verification for Distributed Parameter Systems Using Barrier Functionals\tnoteref{mytitlenote}}
\tnotetext[mytitlenote]{M. Ahmadi was supported by a Clarendon Scholarship and the Sloane-Robinson Scholarship. A. Papachristodoulou is supported in part by EPSRC projects EP/J012041/1, EP/M002454/1 and EP/J010537/1. A preliminary version of this paper was presented at the 2015 American Control Conference, July 1-3, Chicago, IL, USA. Corresponding author A. Papachristodoulou. Tel. +44 1865 2 83036. Fax +44 1865 273010.}

\author{Mohamadreza Ahmadi}
\address{Institute for Computational Engineering and Sciences (ICES), 
University of Texas at Austin, \\Peter O'Donnel, Jr. Building,  201 E 24th St, Austin, TX 78712, USA.}

\author{Giorgio Valmorbida}
\address{Laboratoire des Signaux et Syst\`emes, CentraleSup\'elec, CNRS, Univ. Paris-Sud, Universit\'e Paris-Saclay, 3 Rue Joliot-Curie, Gif sur Yvette 91192, France.}

\author{Antonis Papachristodoulou\corref{mycorrespondingauthor}}
\cortext[mycorrespondingauthor]{Corresponding author}
\ead{antonis@eng.ox.ac.uk}
\address{Department of Engineering Science, University of Oxford, Parks Road, Oxford, OX1 3PJ, UK.}




\begin{abstract}
We study the safety verification problem for a class of distributed parameter systems described by partial differential equations (PDEs), \textit{i.e.,} the problem of checking whether the solutions of the PDE satisfy a set of constraints at a particular point in time. The proposed method is based on an extension of barrier certificates to infinite-dimensional systems. In this respect, we introduce \textit{barrier functionals}, which are functionals of the dependent and independent variables. Given a set of initial conditions and an unsafe set, we demonstrate that if such a functional exists satisfying two (integral) inequalities, then the solutions of the system do not enter the unsafe set. Therefore, the proposed method does not require finite-dimensional approximations of the distributed parameter system. Furthermore, for PDEs with polynomial data, we solve the associated integral inequalities using semi-definite programming (SDP) based on a method that relies on a quadratic representation of the integrands of integral inequalities. The proposed method is illustrated through examples.
\end{abstract}

\begin{keyword}
Safety Verification, Barrier Certificates, Sum-of-Squares Programming, Distributed Parameter Systems
\end{keyword}

\end{frontmatter}


\section{Introduction}

Many real-world engineering systems are described by partial differential equation (PDE) models, which include derivatives with respect to both space and time. For example, mechanics of fluid flows ~\cite{DG95}, dynamics of spatially inhomogeneous robot swarms ~\cite{BermanKN11}, satellite docking systems~\cite{7798766} and the magnetic flux profile in a tokamak ~\cite{6426638} are all described by PDEs. However, compared to systems described by ordinary differential equations (ODEs), the analysis of PDE systems is more challenging. For instance, the solutions to PDEs belong to infinite dimensional (function) spaces, where the norms are not equivalent, as opposed to Euclidean spaces for ODEs. Hence, properties such as stability ~\cite{PP06} and input-output gains ~\cite{Ahmadi2016163} may differ from one norm to another.

One interesting and unresolved problem in the analysis of PDEs is safety verification. That is, given the set of initial conditions, check whether the solutions of the PDE satisfy a set of constraints, or, in other words, whether they are safe with respect to an unsafe set. Reliable safety verification methods are fundamental for designing safety critical systems, such as life support systems ~\cite{4343968},  and wind turbines ~\cite{WSPS13}.
The safety verification problem is well-studied for ODE systems (see the survey paper ~\cite{GLZN09}). Methods based on the approximation of the reachable sets are considered in ~\cite{Kurzhanski2000201} for linear systems and in ~\cite{TMBO03} for nonlinear systems. Another method for safety verification, which does not require the approximation of reachable sets, uses barrier certificates. Barrier certificates ~\cite{P06} were introduced for model invalidation of ODEs with polynomial vector fields and have been used to address safety verification of nonlinear and hybrid systems ~\cite{PJP07} and safety analysis of time-delay systems ~\cite{PJ05}. 
Exponential barrier functions were proposed in ~\cite{steinhardt2012finite} for finite-time regional verification of stochastic nonlinear systems. Moreover, compositional barrier certificates and converse results were studied in ~\cite{SWP12} and \cite{PR05,7236867}, respectively. 

The application of barrier certificates goes beyond just  analysis. Inspired by the notion of control Lyapunov functions~\cite{ARTSTEIN19831163} and Sontag's formula ~\cite{Sontag1989}, Weiland and Allg\"ower~\cite{WA07} introduced control barrier functions (CBFs) and formulated a controller synthesis method that ensures safety with respect to an unsafe set. This has sparked several subsequent studies on control barrier functions ~\cite{7040372, Romdlony201639}.

In this paper, we study the safety verification problem for PDEs using barrier certificates. The proposed method employs a functional of the dependent and independent variables called the barrier functional. We show that the  safety verification problem can be cast as the existence of a barrier functional satisfying a set of integral inequalities. For PDEs with polynomial data, we demonstrate that the associated integral inequalities can be solved using semi-definite programming (SDP) based on the results in ~\cite{valmorbida2016stability}, which were also used in ~\cite{Ahmadi2016163} to solve dissipation inequalities for PDEs and in ~\cite{ahmadi2015convex} for input-output analysis of fluid flows. In this respect, we formulate an S-procedure-like scheme for checking integral inequalities subject to a set of integral constraints. The proposed method is illustrated by two examples.

A preliminary application of the proposed method to bounding nonlinear output functionals of nonlinear time-dependent PDEs was discussed in ~\cite{ahmadi2015barrier}. In this regard, an scheme for bounding linear output functionals of linear stationary PDEs using SDPs was presented in ~\cite{Bertsimas2006} based on moment relaxation techniques. In addition, a moment-relaxation-based method was formulated in ~\cite{Mevissen201110887} to find smooth approximations of the solutions to nonlinear stationary PDEs using a finite-difference discretization of the domain and maximum entropy estimation.

This paper is organized as follows. In the next section, we present some preliminary definitions. In Section 3, we describe a method based on barrier functionals for safety verification of PDEs. In Section 4, we discuss the computational formulation of the barrier functionals method and describe an scheme for verifying integral inequalities subject to integral constraints. We illustrate the proposed results using two examples in Section 5 and conclude the paper in Section 6.

\textbf{Notation:} The $n$-dimensional Euclidean space is denoted by $\mathbb{R}^n$ and the
set of nonnegative reals by $\mathbb{R}_{\ge0}$. The $n$-dimensional set of positive integers is
denoted by $\mathbb{N}^n$, and the $n$-dimensional space of non-negative integers is denoted
by $\mathbb{N}^n_{\ge0}$. We use $M^\prime$ to denote the transpose of matrix $M$. The set of real
symmetric matrices is denoted $\mathbb{S}^n = \{A \in \mathbb{R}^{n\times n} \mid A = A^\prime \}$. The ring of
polynomials on a real variable $x$ is denoted $\mathcal{R}[x]$, and, for $f \in \mathcal{R}[x]$, $deg(f)$
denotes the degree of $f$ in $x$. A domain $\Omega$ is an open subset of $\mathbb{R}^n$ and the
boundary of $\Omega$ is denoted $\partial \Omega$. The space of $k$-times continuous differentiable
functions defined on $\Omega$ is denoted by $\mathcal{C}^k(\Omega)$ and the space of $\mathcal{C}^k(\Omega)$ functions
mapping to a set $\Gamma$ is denoted $\mathcal{C}^k(\Omega;\Gamma)$. For a multivariable function $f(x,y)$,
we use $f(x,\cdot) \in \mathcal{C}^k[x]$ to denote the k-times continuous differentiability of $f$
with respect to variable $x$. If $p \in \mathcal{C}^1(\Omega)$, then $\partial_x p$ denotes the derivative of $p$
with respect to variable $x \in \Omega$. In addition, we adopt Schwartz's multi-index
notation. For $u \in \mathcal{C}^\alpha (\Omega;\mathbb{R}^m)$, $\Omega \in  \mathbb{R}^n$, $\alpha \in \mathbb{N}_{\ge 0}$, defining matrix $A \in \mathbb{N}^{\sigma(m,\alpha) \times n}_{\ge0}$, $\sigma(n,\alpha) = \frac{(n+\alpha)!}{n! \alpha!}$ (denote its $i$th row $A_i$) which contains a set of ordered elements satisfying  $\Sigma_j A_{ij} \le  \alpha$, we have
$$
D^\alpha u := \left(u_1, \partial_x u_1, \ldots, \partial_x^{A_\sigma}u_1, \ldots, u_m, \partial_x u_m, \ldots, \partial_x^{A_\sigma} u_m    \right),
$$
where $\partial_x^{A_i}(\cdot)=\partial_x^{A_{i1}}(\cdot) \cdots \partial_x^{A_{in}}(\cdot)$. We use the same multi-index notation to denote a vector of monomials up to degree $\alpha$ on a variable $x$ as $\eta^\alpha(x)$. For instance,
for $x \in \mathbb{R}^2$, $\eta^2(x) = (1, x_1, x_2, x^2_1, x_1x_2, x_2^2)$. The Hilbert space of functions defined  over the domain $\Omega$ with the norm 
$\|u\|_{\mathcal{W}_\Omega^p} = \left(\int_\Omega \sum_{i=0}^p (\partial_{x_i} u)^\prime (\partial_{x_i} u) ~~\mathrm{d}x\right)^{\frac{1}{2}}$
 is denoted $\mathcal{W}_\Omega^p$. By $f \in \mathcal{L}^2(\Omega;\Gamma)$, we denote a square integrable function mapping $\Omega \subseteq \mathbb{R}^n$ to $\Gamma \subseteq \mathbb{R}^m$. Also, for an operator $\mathscr{A}$, $Dom(\mathscr{A})$ and $Ran(\mathscr{A})$ denote its domain and range, respectively. The notation $\lceil \cdot \rceil$ denotes the ceiling function.

\section{Preliminaries} \label{sec:safety}

In this section, we present some definitions and preliminary results. We study a class of forward-in-time PDE systems. Let $ \mathcal{U}$ be a Hilbert space. Consider the following differential equation
\begin{equation}  
\label{eq:pde}
\begin{cases}
\partial_t u(t,x)= \mathscr{F} u(t,x), \quad x\in \Omega \subset\real^n,~t \in [0,T],\\
y(t) = \mathscr{H} u(t,x)\\
u(0,x)=u_0(x) \in \mathcal{U}_0\subset Dom(\mathscr{F})\\
u\in \mathcal{U}_b  
\end{cases}
\end{equation}
where $\mathcal{U}_b$ is a subspace of $\mathcal{U}$, the state-space of system~\eqref{eq:pde}, defined by the boundary conditions, $\mathscr{H}: \mathcal{U} \rightarrow \real$ and $Dom(\mathscr{H})\supseteq \mathcal{U}$, the state-space of system~\eqref{eq:pde}. It is assumed that \eqref{eq:pde} is well-posed.  \ref{app:wellposed} reviews some aspects of the well-posedness of PDEs. While these results are important, studying the well-posedness of system~\eqref{eq:pde} is beyond the scope of the current paper.

 We call the set  
$$
 \mathcal{Y}_u = \big \{ u \in  \mathcal{U} \mid \mathscr{H}  u \le 0  \big\},
$$
the \emph{unsafe set}. 


Consider the following properties of trajectories related to an initial set $\mathcal{U}_0$ and an unsafe set $\mathcal{Y}_u$.

\begin{defi} [Safety  at Time $T$] \label{def:safe}
Let $u\in~\mathcal{U}$. For a set $ \mathcal{U}_0 \subseteq  \mathcal{U}$, an unsafe set $ \mathcal{Y}_u$, satisfying $ \mathcal{U}_0 \cap \mathcal{Y}_u = \emptyset$, and a positive scalar $T$, system~\eqref{eq:pde}  is \emph{$ \mathcal{Y}_u$-safe at time $T$}, if the solutions $u(t,x)$ of system~\eqref{eq:pde} satisfy $y(T) \notin  \mathcal{Y}_u$ for all $u(0,x) \in \mathcal{U}_0$.
\end{defi}

\begin{defi} [Safety ] \label{def:safe2}
System~\eqref{eq:pde} is  \emph{$\mathcal{Y}_u$-safe}, if it is safe with respect to $\mathcal{Y}_u$ in the sense of Definition~\ref{def:safe} for all $T>0$.
\end{defi}

We are interested in solving the following problem:

\begin{problem}\label{safetyproblem}
Given sets $\mathcal{Y}_u$, $ \mathcal{U}_0$  and a constant $T>0$, verify that system~\eqref{eq:pde} is  $ \mathcal{Y}_u$-safe at time $T$.
\end{problem}

To this end, we define a time-dependent functional of the states of the PDE and time 
\begin{equation} \label{eq:barrierfunctional}
B(t,u) = \mathscr{B}(t)u,
\end{equation}
 where $\mathscr{B}(t) : Dom(\mathscr{B}) \to \mathbb{R}$. We refer to this functional as the {barrier functional}.  Note that this  extension of barrier certificates~\cite{P06} enables us to address sets that are defined on  infinite-dimensional spaces.   In the subsequent section, we show that the  barrier functional provides the means to characterize a barrier between the set of initial conditions and the unsafe set.

\section{Barrier Functionals for Safety Verification of PDEs}

In this section, we present conditions to obtain certificates that trajectories starting in the set $\mathcal{U}_0$ are $\mathcal{Y}_u$-safe at a particular time instant $T$. Such a formulation also allows obtaining performance estimates
whenever the unsafe set represents a performance index.
 
Next, we provide a solution to Problem~\ref{safetyproblem}  based on the construction of barrier functionals satisfying a set of inequalities. 
 
\begin{thm}[Safety Verification for Forward PDE Systems] \label{thm:fwd}
Consider the PDE system described by~\eqref{eq:pde}. Let $u \in   \mathcal{U}_b$. Given  a set of initial conditions $ \mathcal{U}_0 \subseteq  \mathcal{U}_b$, an unsafe set~$ \mathcal{Y}_u$, such that $ \mathcal{U}_0 \cap  \mathcal{Y}_u = \emptyset$,  and a constant $T>0$, if there exists a barrier functional $B(t,u(t,x)) \in \mathcal{C}^1[t]$ as in \eqref{eq:barrierfunctional},
 such that the following inequalities hold
 \begin{subequations}
\label{eq:confwd}
\begin{equation} \label{con1}
B(T, u(T,x)) - B(0, u_0(x)) >0,\quad \forall u(T,x) \in  \mathcal{Y}_u,~\forall u_0 \in  \mathcal{U}_0,
\end{equation}
\begin{equation} \label{con2}
 \frac{~\mathrm{d}B(t,u(t,x))}{~\mathrm{d}t} \le 0,\quad \forall t\in[0,T],~\forall u \in  \mathcal{U}_b,
\end{equation}
\end{subequations}
where $\frac{~\mathrm{d}(\cdot)}{~\mathrm{d}t}$ denotes the total derivative, along the solutions of~\eqref{eq:pde}, then the solutions of~\eqref{eq:pde} are $ \mathcal{Y}_u$-safe  at time $T$ (cf. Definition~\ref{def:safe}).
\end{thm}
\begin{pf}
The proof is by contradiction. Assume there exists a solution of \eqref{eq:pde} such that, at time $T$, $u(T,x) \in  \mathcal{Y}_u$ and inequality~\eqref{con1} holds.  From \eqref{con2}, it follows that
\begin{equation} \label{esasq1}
 \frac{~\mathrm{d}  B(t,u(t,x))}{~\mathrm{d}t}    \le 0,
\end{equation}
for all  $t\in[0,T]$, and  $u \in  \mathcal{U}$. Integrating both sides of~\eqref{esasq1} with respect to $t$ from $0$ to $T$ yields
\begin{equation*} 
\int_0^{T} \frac{~\mathrm{d}  B(t,u)}{~\mathrm{d}t} ~\mathrm{d}t =  B(T, u(T,x)) - B(0, u(0,x)) \le  0.
\end{equation*}
for all $u \in  \mathcal{U}$. This  contradicts \eqref{con1}.
\end{pf}

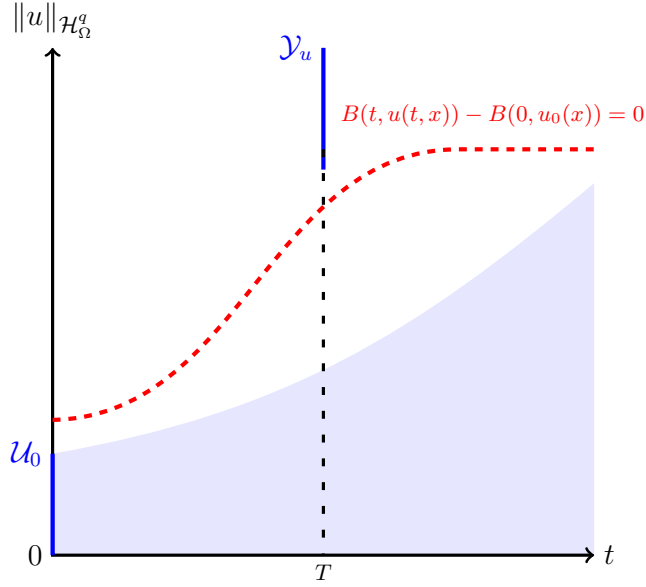
\begin{figure}
  \centering
\begin{tikzpicture}[scale=0.9, every node/.style={scale=0.95}]
\path[fill=blue!10] (0,0)--(0,1.5) to [out=10,in=220] (8,5.5) --(8,0);
\draw[very thick,<->] (8,0) node[font=\fontsize{13}{58}\sffamily\bfseries,right]{$t$} -- (0,0) node[font=\fontsize{13}{58}\sffamily\bfseries,left ]{$0$} -- (0,7.5) node[font=\fontsize{13}{58}\sffamily\bfseries,above]{$\| u \|_{ \mathcal{H}^q_{\Omega}}$};
\draw[dashed,red, ultra thick] (0,2) to [out=1,in=180] (6,6) --(8,6);
\draw[ultra thick,blue] (4,5.7)--(4,7.5) node[font=\fontsize{12}{58}\sffamily\bfseries,left ]{$ \mathcal{Y}_u$};
\draw[very thick,loosely dashed] (4,6)--(4,0) node[font=\fontsize{9}{58}\sffamily\bfseries,below ]{$T$};
\draw[ultra thick, blue] (0,0) --(0,1.5)node[font=\fontsize{12}{58}\sffamily\bfseries,left ]{$ \mathcal{U}_0$};
\draw [red] (6.5,6.8) node[font=\fontsize{9}{58}\sffamily\bfseries,below ]{$B(t,u(t,x))-B(0,u_0(x))=0$};
\end{tikzpicture}
  \caption{Illustration of a barrier functional for a  PDE system: any solution $u(t,x)$ with  $u(0,x) \in \mathcal{U}_0$ (depicted by the shaded area) satisfies $u(T,x) \notin  \mathcal{Y}_u$. The system is $ \mathcal{Y}_{u}$-safe at time $t=T$ but not for  $\forall t>0$.}\label{tikzfig1}
\end{figure}

\begin{remark} The level sets of~$B(t,u(t,x))-B(0, u_0(x))$  represent barrier surfaces in the $ \mathcal{U}$ space separating $ \mathcal{U}_0$ and $ \mathcal{Y}_u$ such that no solution of \eqref{eq:pde}  starting from $ \mathcal{U}_0$ is in $ \mathcal{Y}_u$ at time $T$ (hence, the term ``barrier functional''). This property is illustrated in Figure~\ref{tikzfig1}.
\end{remark}

Theorem~\ref{thm:fwd} is concerned with conditions for safety verification with respect to the unsafe set $ \mathcal{Y}_u$ at a particular time $T>0$. The next corollary follows from Theorem~\ref{thm:fwd} and gives conditions for safety verification with respect to an unsafe set $ \mathcal{Y}_u$ for all time $t>0$. In this case, the barrier functional can be independent of $t$.
\begin{cor} \label{cor:fwd}
Consider the PDE system described by~\eqref{eq:pde}.  Assume  $u \in   \mathcal{U}_b$. Given an unsafe set~$\mathcal{Y}_u \subset  \mathcal{U}$, such that~$ \mathcal{U}_0 \cap  \mathcal{Y}_u = \emptyset$, if there exists a barrier functional $B(u(t,x))$ as in \eqref{eq:barrierfunctional}
 such that
 \begin{subequations}\label{eq:cor:fwdcon}
 \begin{equation} \label{eq:cor:fwdcon1}
B(u(t,x)) - {B}(u_0(x)) > 0,\quad \forall u \in  \mathcal{Y}_u,~ \forall u_0 \in  \mathcal{U}_0,
 \end{equation}
 \begin{equation} \label{eq:cor:fwdcon2}
 \frac{~\mathrm{d}B(u(t,x))}{~\mathrm{d}t} \le 0,\quad \forall u \in  \mathcal{U},
 \end{equation}
 \end{subequations}
 along the solutions of \eqref{eq:pde}, then the solutions of  PDE~\eqref{eq:pde} are $ \mathcal{Y}_u$-safe (cf. Definition~\ref{def:safe2}).
\end{cor}

\begin{pf}
The proof follows the same lines as the proof of Theorem~\ref{thm:fwd}. Assume that there exists a solution $u(t,x)$ to~\eqref{eq:pde} such that, for some $t>0$, we have $u(t,x) \in \mathcal{Y}_u$. Then, from \eqref{eq:cor:fwdcon1}, it follows that $B(u(t,x)) - {B}(u_0(x)) > 0$. On the other hand, integrating inequality \eqref{eq:cor:fwdcon2} from $0$ to $t$ implies that $B(u(t,x)) - {B}(u_0(x)) \le 0$, which is a contradiction. Thus, since $t$ is arbitrary, the solutions to \eqref{eq:pde} are $\mathcal{Y}_u$-safe for all time.
\end{pf}

We conclude this section by illustrating Corollary~\ref{cor:fwd} with an analytical example that uses a barrier functional to bound a performance index. 
\begin{exam}\textbf{(Performance Bounds)} Consider the heat equation defined over a domain \mbox{$\Omega \subset \mathbb{R}^2$} with smooth boundary
\begin{equation} \label{eq:exampleheat}
\partial_t u = \Delta u,~~~~x \in \Omega,~t>0,
\end{equation}
subject to boundary conditions $u |_{\partial \Omega} =0$ and  
\begin{equation} \label{xcxz}
u(0,x) \in  \mathcal{U}_0 = \left \{ u_0 \in \mathcal{U}  \mid \int_{\Omega} |\nabla u_0 |^2 ~\mathrm{d}\Omega \le 1 \right\}.
\end{equation}
where $\Delta$ is the Laplacian operator. The output mapping is given by
$$
y(t) =   \gamma^2 - \int_\Omega u^2(t,x)~\mathrm{d}\Omega,
$$
where $\gamma\geq0$. Then, the unsafe set is described as
$
 \mathcal{Y}_u =\left \{ u \in \mathcal{U}  \mid  y(t) = \gamma^2 - \int_\Omega u^2(t,x)~\mathrm{d}\Omega < 0 \right \}
$. We are interested in finding the minimum $\gamma$ such that no solution of~\eqref{eq:exampleheat}  enters $ \mathcal{Y}_u$ for all $u(0,x) \in  \mathcal{U}_0$.

We consider the  barrier functional~\eqref{eq:barrierfunctional} with
$$
\begin{array}{cc}
\mathscr{B}: & \mathcal{W}^1_\Omega \rightarrow \real_{\geq 0} \\ & u \mapsto \int_\Omega (\nabla u)^\prime  \nabla u  ~\mathrm{d}\Omega,
\end{array}
$$
that is, $B(u(t,x))= \int_\Omega (\nabla u)^\prime  \nabla u  ~\mathrm{d}\Omega$. We first check inequality \eqref{eq:cor:fwdcon2} along the solutions of \eqref{eq:exampleheat}:
\begin{align*} 
 \frac{~\mathrm{d} B(u(t,x))}{~\mathrm{d}t}  &= \int_\Omega 2 \nabla u \partial_t \left( \nabla u  \right) ~\mathrm{d}\Omega  = 2 \left(\nabla u \partial_t   u \right)|_{\partial \Omega} - 2 \int_\Omega \Delta u \partial_t   u  ~\mathrm{d}\Omega \\ &= -2 \int_\Omega  \left( \Delta u \right)^2  ~\mathrm{d}\Omega \le 0,
\end{align*}
where, in the second equality above, integration by parts and, in the third equality, the boundary conditions are used. Thus, inequality \eqref{eq:cor:fwdcon2} is satisfied. At this point, let us check inequality~\eqref{eq:cor:fwdcon1}. We have\begin{eqnarray} \nonumber
B(u(t,x))-B(u_0) &=&
\int_\Omega |\nabla u |^2 ~\mathrm{d}\Omega - \int_\Omega |\nabla u_0 |^2 ~\mathrm{d}\Omega \nonumber  \ge \int_\Omega |\nabla u  |^2  ~\mathrm{d}\Omega -1 \nonumber \\
 &\ge& C(\Omega) \int_\Omega  u^2~\mathrm{d}\Omega -1, \nonumber
\end{eqnarray}
where $u_0 \in  \mathcal{U}_0$ as in~\eqref{xcxz} is applied to obtain the first inequality, $C(\Omega)>0$, and the Poincar\'{e} inequality~~\cite{PW60} is used in the second inequality. Then, it follows that whenever $\gamma^2 > \frac{1}{C(\Omega)}$, we have
$
B(u(t,x))-B(u_0)>0,
$ 
and thus, from Theorem~\ref{thm:fwd}, system~\eqref{eq:exampleheat} is $ \mathcal{Y}_u$-safe. Therefore, it holds that $y \notin  \mathcal{Y}_u$, which implies $y(t) = \gamma_{min}^2 - \int_\Omega u^2~\mathrm{d}\Omega \ge 0$, \textit{i.e.,} $
\gamma_{min}^2 \ge \int_\Omega u^2~\mathrm{d}\Omega,$ 
where $\gamma_{min}^2 = \frac{1}{C(\Omega)}$. For example, whenever $\Omega = \{ (x,y)\in \mathbb{R}^2 \mid |x+y|<1 \}$, we obtain $\gamma^2 =\frac{{2}}{\pi^2}$.
\end{exam}

{
Note that  in the  example above, for which the barrier functional is fixed, the verification of the inequalities in Theorem 3.1 can be involved. Moreover, when the dynamics are nonlinear or spatially varying and for more general  initial/unsafe  sets, the selection of the barrier certificate candidate is less obvious. Therefore, we bring forward a numerical method to automate the verification steps and hence the construction of the barrier certificates in the next section. }

\section{Construction of Barriers Functionals} \label{sec:computational}

In this section, we study a specific class of barrier functionals and  particular sets $\mathcal{U}_0$ and $\mathcal{Y}_u$,  for which the inequalities~\eqref{eq:confwd} become integral inequalities. For the case of polynomial data, the verification of the inequalities can be cast as constraints of an SDP. Furthermore, we set $\Omega = (0, 1)$. Note that any bounded open subset of the real line can be mapped into this domain\footnote{A domain $(a,b)$ can be mapped to $(0,1)$ by the following change of variables $$\bar{x}=\frac{x-a}{b-a}.$$}.

In the previous sections, the barrier functionals were only assumed to be continuously differentiable with respect to time. In order to present computational tools based on SDPs, hereafter we assume that the barrier functional takes the form of an integral functional.

\section{Verifying Integral Inequalities with Integral Constraints}

In order to check inequalities~\eqref{eq:confwd} and~\eqref{eq:cor:fwdcon} based on the method proposed in ~\cite{valmorbida2016stability}, we require verifying an integral inequality subject to a number of integral constraints. That is, we need to solve the following class of problems

\begin{eqnarray}\label{optmiconsteq}
& \int_0^1 f_i(t,x,D^\alpha u) ~\mathrm{d}x\ge 0,&\nonumber \\
&\text{subject to}&\nonumber \\
&\int_0^1 s_i(t,x,D^\alpha u)~\mathrm{d}x \ge 0,~i=1,2,\ldots,r.&
\end{eqnarray}
where $u : \mathbb{R}_{\ge0} \times \Omega \to  \mathbb{R}^n$, $f_i, s_i \in \mathcal{R}[t, x, D^\alpha u]$ and $max(deg(s_i), deg(f_i)) = k$. Let $\sigma(n,k):=\frac{(n+k-1)!}{(n-1)!k!}$. Then, we can represent $f_i$ and $s_i$ as the following quadratic-like forms
$$
f_i(t,x,D^\alpha u) = \left(\eta^{\left\lceil \frac{k}{2} \right\rceil} (D^\alpha u)\right) ^\prime  F_i(t,x)~\eta^{\left \lceil \frac{k}{2} \right \rceil} (D^\alpha u)
$$ 
$$
s_i(t,x,D^\alpha u) = \left(\eta^{\left\lceil \frac{k}{2} \right\rceil} (D^\alpha u)\right) ^\prime  S_i(t,x)~\eta^{\left \lceil \frac{k}{2} \right \rceil} (D^\alpha u)
$$
with $F_i,S_i: \mathbb{R}_{\ge0}\times \Omega \rightarrow \mathbb{S}^{\sigma(n\alpha, \left\lceil \frac{k}{2} \right \rceil) }$.

The approach we develop here is reminiscent of S-procedure ~\cite{doi101137S003614450444614X} for LMIs. The S-procedure provides conditions under which a particular quadratic inequality holds subject to some other quadratic inequalities (for example, within the intersection of several ellipsoids). Similar conditions for checking polynomial inequalities within a semi-algebraic set were developed in~\cite{PP05, Par00} thanks to Putinar's Positivstellensatz~\cite[Theorem 2.14]{Las09}. However, current machinery for including integral constraints includes multiplying the integral constraint and subtracting it from the inequality (see Proposition 9 in ~\cite{PP05}). In the following, we propose an alternative to the latter method that can be used to verify the feasibility problem (8).

Consider the following set of integral constraints
\begin{equation}\label{p9303od03}
\mathcal{S} = \left\{ u \in \mathcal{C}^{\alpha} (\Omega;\mathbb{R}^n) \mid \int_\Omega s_i(t,x,D^\alpha u)~\mathrm{d}x \ge 0,~i=1,2,\ldots,r \right\}.
\end{equation}
Note that in this setting, we can also represent sets as $\left\{ u  \mid \int_\Omega g(t,x,D^\alpha u)~\mathrm{d}x =0\right\}$ by selecting $s_1 = g$ and $s_2 =-g$.

Define 
\begin{equation} \label{etrerter}
v_i(t,x) := \int_0^x s_i(t,x,D^\alpha u)~\mathrm{d}x,
\end{equation}
satisfying 
\begin{equation} \label{fsdsdffccc}
\begin{cases}
v_i(t,0) = 0,\\
\partial_x v_i(t,x) - s_i(t,x,D^\alpha u(t,x))=0,
\end{cases}
\end{equation}
for $i=1,2,\ldots,r$. Using~\eqref{etrerter}, we can represent $\mathcal{S}$ as 
$$
\mathcal{S} = \left\{ u \in  \mathcal{C}^{\alpha}_\Omega \mid v_i(t,1) \ge 0,~i=1,2,\ldots,r  \right\}.
$$
\begin{lem}\label{lem:integralcons22233}
Consider problem~\eqref{optmiconsteq} and let $t \in \mathcal{T} \subseteq \mathbb{R}_{\ge0}$. Let $v(t,x) = \left[\begin{smallmatrix} v_1(t,x) &\cdots& v_r(t,x)  \end{smallmatrix}\right]^\prime$ and $s(x,D^\alpha u) = \left[\begin{smallmatrix} s_1(x,D^\alpha u) & \cdots & s_r(x,D^\alpha u)\end{smallmatrix} \right]^\prime$. If there exists a vector function $m: \mathcal{T}  \times \Omega \to \mathbb{R}^r$ and a vector $n \in \mathbb{R}^r_{\ge 0}$ such that
\begin{equation} \label{xssddscx}
\int_0^1 f_i(t,x,D^\alpha u)~\mathrm{d}x -n^\prime v(t,1) + \int_0^1 m^\prime(t,x)\Big(\partial_x v(t,x) - s\left(t,x,D^\alpha u(t,x)\right) \Big)~\mathrm{d}x >0,
\end{equation}
for all  $u \in \mathcal{U}$ and all $t \in \mathcal{T} $, then~\eqref{optmiconsteq} is satisfied. 
\end{lem}
\begin{pf}
From~\eqref{fsdsdffccc}, we have that for any $m: \mathbb{R}_{\ge0} \times \Omega \to \mathbb{R}^r$
$$
m^\prime(t,x)\left(\partial_x v(t,x) - s(t,x,D^\alpha u(t,x)) \right)=0,\quad \forall x \in \Omega.
$$
Hence, since $v$ and $u$ are related according to~\eqref{fsdsdffccc}, we obtain
$$
\int_0^1 m^\prime(t,x)\left(\partial_x v(t,x) - s(t,x,D^\alpha u(t,x)) \right)~\mathrm{d}x=0.
$$
Consequently, if inequality~\eqref{xssddscx} is satisfied, we infer 
$$
\int_0^1 f_i(t,x,D^\alpha u)~\mathrm{d}x> n^\prime v(t,1),\quad \forall t \in \mathcal{T}.
$$
Finally, since $n^\prime v(t,1)\ge 0$, for all $u \in \mathcal{S}$, we conclude that ~\eqref{optmiconsteq} holds. 
\end{pf}

Note that inequality (12) can be checked using the method discussed in ~\cite{valmorbida2016stability}.
In order to incorporate the integral constraints, we introduced the (dummy)
dependent variables $v_i(t,x)$, satisfying (11), and their partial derivative with
respect to $x$. 

\subsection{Computational Formulation}

We impose the following structure for the barrier functionals\begin{equation}
B(t,u) = \int_\Omega \left(\eta^{\left\lceil \frac{k}{2} \right\rceil} (D^\alpha u)\right)^\prime~ \bar{B}(t,\theta)~\left(\eta^{\left\lceil \frac{k}{2} \right\rceil} (D^\alpha u)\right) ~\mathrm{d} \theta
\label{eq:barrierquad}
\end{equation}
where $\Omega = (0,1)$,  $\bar{B}: \real_{\geq 0}\times \Omega \rightarrow \real^{\sigma(n\alpha,{\left\lceil \frac{k}{2} \right\rceil}) \times \sigma(n\alpha,{\left\lceil \frac{k}{2} \right\rceil})}$, $\bar{B}(t,x)\in \mathcal{C}^{1}[t], \forall x\in\Omega$, and the following quadratic-like structures for the unsafe and the initial sets
\begin{subequations}
\label{eq:quadsets}
\begin{equation}\label{eq:Yuquad}\mathcal{Y}_u =  \left\lbrace u \in \mathcal{U} \mid   \int_\Omega \left(\eta^{\left\lceil \frac{k}{2} \right\rceil} (D^\alpha u)\right)^\prime~ Y(t,\theta) ~\left(\eta^{\left\lceil \frac{k}{2} \right\rceil} (D^\alpha u)\right)   ~\mathrm{d}\theta \geq 0  \right\rbrace,\end{equation} and the  set of initial conditions 
\begin{equation}\mathcal{U}_0 =  \left\lbrace u_0 \in \mathcal{U} \mid  \int_\Omega \left(\eta^{\left\lceil \frac{k}{2} \right\rceil} (D^\alpha u)\right)^\prime~U_0(t,\theta) ~\left(\eta^{\left\lceil \frac{k}{2} \right\rceil} (D^\alpha u)\right)  ~\mathrm{d}\theta \geq 0 \right\rbrace. \label{eq:U0quad}\end{equation}
\end{subequations}
where, $Y: \real_{\geq 0}\times \Omega \rightarrow \real^{\sigma(n\alpha,{\left\lceil \frac{k}{2} \right\rceil}) \times \sigma(n\alpha,{\left\lceil \frac{k}{2} \right\rceil})}$ and  $U_0: \real_{\geq 0}\times \Omega \rightarrow \real^{\sigma(n\alpha,{\left\lceil \frac{k}{2} \right\rceil}) \times \sigma(n\alpha,{\left\lceil \frac{k}{2} \right\rceil})}$.

The following proposition applies   Lemma~\ref{lem:integralcons22233} to formulate integral inequalities to verify the conditions of Theorem~\ref{thm:fwd}  considering  barrier functional~\eqref{eq:barrierquad}. In this case, the constraint set $\mathcal{S}$ as defined in \eqref{p9303od03} is given by  $\mathcal{S} = \mathcal{Y}_u \cup \mathcal{U}_0$, with the sets  in~\eqref{eq:quadsets} defined as
\begin{align}
\label{eq:gprop}
s_1(t,x,D^\alpha u) &= \left(\eta^{\left\lceil \frac{k}{2} \right\rceil} (D^\alpha u)\right)^\prime Y(t,x) \left(\eta^{\left\lceil \frac{k}{2} \right\rceil} (D^\alpha u)\right),\nonumber \\
s_2(t,x,D^\alpha u) &= \left(\eta^{\left\lceil \frac{k}{2} \right\rceil} (D^\alpha u)\right)^\prime U_0(t,x) \left(\eta^{\left\lceil \frac{k}{2} \right\rceil} (D^\alpha u)\right). 
\end{align}

\begin{prop}
\label{prop:intineqbarriers}
If there exist $\bar{B}: [0,T] \times \Omega \rightarrow \real^{\sigma(n\alpha,{\left\lceil \frac{k}{2} \right\rceil}) \times \sigma(n\alpha,{\left\lceil \frac{k}{2} \right\rceil})}$ or $B(t,u)$ as in~\eqref{eq:barrierquad}, $m: \mathcal{T} \times \Omega \rightarrow \real^2$ and $n \in \real^{2}_{\geq0}$  such that the following inequalities are satisfied
\begin{subequations}
\begin{multline}
B(T,u(T,x))  - B(0,u_0)    - n^\prime \left[\begin{smallmatrix} v_2(T,1) \\ v_1(0,1)  \end{smallmatrix}\right]  \\
+ \int_{\Omega}^{} \left[\begin{smallmatrix} m_2(T,\theta) \\ m_1(0,\theta)  \end{smallmatrix}\right]^\prime \left[\begin{smallmatrix} \partial_\theta v_2(T,\theta)  -   s_2(T,x,D^\alpha u(T,\theta))\\ \partial_\theta v(0,\theta)  -   s_1(0,x,D^\alpha u_0(\theta))  \end{smallmatrix}\right]  ~\mathrm{d} \theta > 0,
\end{multline}
with $s_1$ and $s_2$ as defined by~\eqref{eq:gprop} and $v_1$ and $v_2$ as defined by~\eqref{etrerter}, and 
\begin{multline}
\label{eq:Bder}
 \int_\Omega \bigg( \left(\eta^{\left\lceil \frac{k}{2} \right\rceil} (D^\alpha u)\right)^\prime  \partial_t \bar{B}(t,\theta) \left(\eta^{\left\lceil \frac{k}{2} \right\rceil} (D^\alpha u)\right) \\+ 2 \left(\eta^{\left\lceil \frac{k}{2} \right\rceil} (D^\alpha u)\right)^\prime   \bar{B}(t,\theta) \nabla \left(\eta^{\left\lceil \frac{k}{2} \right\rceil} (D^\alpha u)\right)^\prime  \partial_t (D^\alpha u)\bigg)~~\mathrm{d}\theta\leq 0,
\end{multline}
\label{eq:barrier_poly}
\end{subequations}
$\forall t \in \left[0,T\right]$, $\forall u \in \mathcal{U}$,  then~\eqref{eq:confwd} holds. 
\end{prop}

A method to solve integral inequalities as (16) was proposed in ~\cite{valmorbida2016stability} (also see ~\cite{valmorbida2015convex} for the formulation for $\Omega \subset \mathbb{R}^2$). In the proposed method, the problem  of checking an integral inequality is cast as the problem of solving a differential linear matrix inequality. Such a formulation is possible thanks to the use of quadratic-like expressions as in (13), (14). Furthermore, it is demonstrated that, for polynomial data, the corresponding differential matrix inequalities can be converted to a Sum-of-Squares (SOS) program, which is then cast as an SDP.  The numerical results presented in the next section consider the problem data to be polynomial, i.e., the functions $\bar{B}$, $m$, $Y$, $U_0$ appearing in the inequalities of Proposition~\ref{prop:intineqbarriers} are polynomials on variables $t$ and $x$, and the operator $\mathscr{F}$ in~\eqref{eq:pde} may be nonlinear and defined by a polynomial on $u$ and its spatial derivatives with coefficients that are polynomials on the spatial variables. The formulation of the SDPs can be automated and a plug-in to SOSTOOLS~\cite{valmorbida2015introducing} has been developed.

\section{Examples} \label{sec:examples}

We now illustrate the proposed results with two numerical examples. The first example is associated with the option pricing problem from quantitative finance. The second example concerns a diffusion-reaction-convection PDE. The numerical results given in this section were obtained using SOSTOOLS v. 3.00 ~\cite{PAVPSP13} and the associated SDPs were solved using SeDuMi v.1.02 ~\cite{Stu98}.


\subsection{Example 1: Option Pricing}


Consider the following linear PDE
\begin{equation} \label{eq:BS}
\partial_t u(t,s) = \frac{\sigma^2 s^2}{2} \partial_s^2 u(t,s) + rs \partial_s u(t,s)-ru(t,s),~(t,s)\in [0,T]\times [0,\bar{s}],
\end{equation}
which is  the (forward) Black-Scholes equation for a non-dividend-paying stock (see~\cite[p. 331]{Hu15}). For the European call option the terminal and the boundary conditions are given as
 \begin{equation*}
 \begin{cases}
 u(T,s)= f(s)=\max\left\{s-K, 0  \right\},\\
 u(t,0)= 0, \\
 u(t,\bar{s})=\bar{s},
 \end{cases}
 \end{equation*}
 where $K>0$ is the strike price. Assuming the stock is at-the-money, $f(s)=s-K$. The parameter values for a European call option~\cite[p. 338]{Hu15} are described as
 \begin{equation*}
 \begin{cases}
 T= 0.5~\text{(years)},~K=\$40,\\
 r=0.1,~\sigma=0.2.
 \end{cases}
 \end{equation*}
 We are interested in checking the safety of the solutions to (17) such that the average option price $\frac{1}{s} \int_0^{\bar{s}} u(T,s)~\mathrm{d}s$ does not exceed some price $\gamma$. Of course, a minimization over $\gamma$ gives us an estimate on the actual average option price. 
 
%
To this end, we define
 $$
 \mathcal{Y}_u = \left\lbrace u \in \mathcal{L}^{1}_{[0,\bar{s}]}  \mid \frac{1}{\bar{s}} \int_0^{\bar{s}} u(0,\theta) ~\mathrm{d}\theta  - \gamma  \ge 0 \right\rbrace.
 $$
 Consider the following barrier functional
 \begin{equation*}
 B(t,u(t,x)) = \int_0^{\bar{s}} b(t,\theta) u^2(t,\theta) ~\mathrm{d}\theta.
 \end{equation*}
where $b \in \mathcal{R}[t,\theta]$. Using Proposition~\ref{prop:intineqbarriers} and a minimization over $\gamma$, we obtain the upper bounds on the average option price as given in Table~\ref{t11kjhkjjkkj}. In these numerical experiments, we set $deg(b) = deg(m)$. The actual upper bound obtained from the solution to (17)~\cite[p. 76]{9780511812545} for the average option price is $18.227$. As it can be observed from the table, increasing the degree of the involved polynomials improves the accuracy of $\gamma^\star$. The constructed barrier
 functional certificate of degree 6 is given in Appendix B.

\begin{table}[!t]
\caption{Bounds on the average option price.}
\label{t11kjhkjjkkj}
\centering
\begin{tabular}{c||c|c|c|c|c|c}

\bfseries $deg(b)$ & 1 & 2 & 3& 4 & 5 & 6  \\
\hline
\bfseries $\gamma^\star$ & 44.4285 & 26.1093 & 22.7489 & 19.5572 & 18.8264 & 18.2391   \\
\end{tabular}
\end{table}

\subsection{ Diffusion-Reaction-Convection PDE}

\begin{figure}[!tb]
\begin{psfrags}
     \psfrag{t}[l][l]{\footnotesize $t$}
     \psfrag{x}[l][l]{\scriptsize $x$}
     \psfrag{u}[c][l]{\footnotesize $u(t,x)$}
\epsfxsize=10cm
\centerline{\epsffile{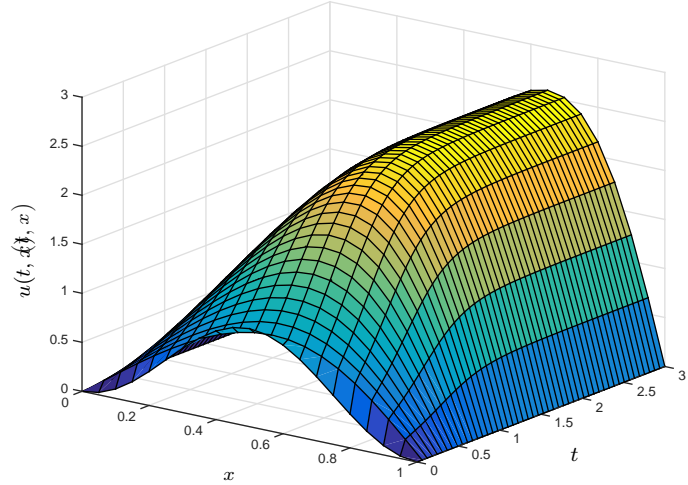}}
\end{psfrags}
\caption{The solution to PDE~\eqref{eq:examplemain} for $\lambda=1.2\pi^2$.\label{fig3}}
\end{figure}

Consider the following nonlinear PDE 
\begin{equation} \label{eq:examplemain}
\partial_t u = \partial_x^2 u + \lambda u - 2 u \partial_x u,~x\in (0,1),~t>0
\end{equation}
where $\lambda>0$, and $u(t,0)=u(t,1)=0$.  Due to the presence of a nonlinear convection term, the solutions with $\lambda \ge \pi^2$ (otherwise unstable) may converge to a different stationary solution. Figure~\ref{fig3} depicts a solution to PDE~\eqref{eq:examplemain} with $\lambda > \pi^2$. {This stems from the fact that the nonlinear convection term transfers low wave number components of the solutions to the high wave number ones for which the diffusion term has a stabilizing effect. This phenomenon appears in the solutions of some important nonlinear PDEs, including Kuramoto-Sivashinsky equation~\cite{FANTUZZI201523}, Burgers Equation~\cite{KRSTIC1999123} and the KdV equation~\cite{doi16M1061837}. }

\begin{figure}[!tb]
\begin{psfrags}
     \psfrag{t}[l][l]{\footnotesize $t$}
     \psfrag{u100000}[l][l]{\scriptsize $u_0 = 2x^2(1-x)^2$}
     \psfrag{u2}[l][l]{\scriptsize $u_0 = 0.4x(e^x - e)$}
     \psfrag{u3}[l][l]{\scriptsize $u_0 = x(1-x)$}
     \psfrag{u4}[l][l]{\scriptsize $u_0 = \ln (\frac{11}{20}x(1-x)+1)$}
     \psfrag{u}[c][l]{\footnotesize $\|u\|_{\mathcal{W}^{1}_{(0,1)}}~~~$}
\epsfxsize=10cm
\centerline{\epsffile{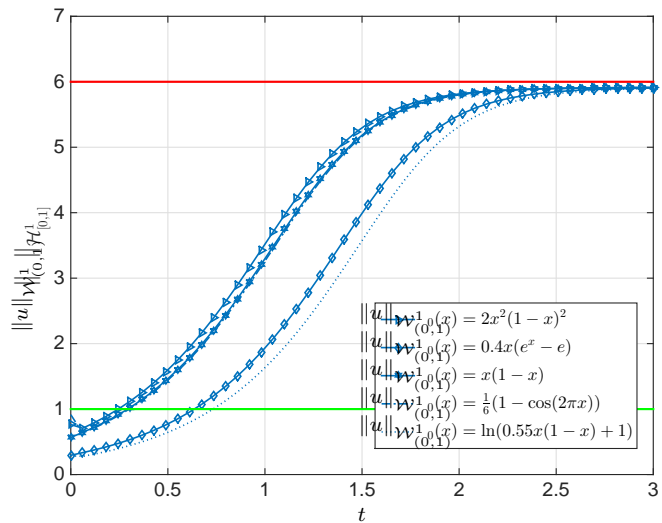}}
\end{psfrags}
\caption{The evolution of $\mathcal{H}^1_{(0,1)}$-norm of solutions to \eqref{eq:examplemain} with  $\lambda=1.2\pi^2$ for  different initial conditions. The red and green lines show the boundaries of $ \mathcal{Y}_{u}$ and $ \mathcal{U}_0$, respectively.\label{fig2}}
\end{figure}

We are interested in computing the maximum value for parameter $\lambda$, such that the solutions starting in
\begin{equation} \label{eq:u00}
 \mathcal{U}_0 = \left\lbrace u_0   \mid \int_0^1 \left( u_0^2 + (\partial_\theta u_0)^2  \right) ~\mathrm{d}\theta \le 1  \right\rbrace,
\end{equation}
which implies $\| u_0 \|_{\mathcal{H}^1_{(0,1)} } \le 1$, do not enter the set
$$ \mathcal{Y}_u = \left\lbrace u   \mid \int_0^1 \left( u^2 + (\partial_\theta u)^2  \right) ~\mathrm{d}\theta \ge (6)^2   \right\rbrace,$$ 
\textit{i.e.,} $\| u \|_{\mathcal{H}^1_{(0,1)}} \ge {6}$ for all $t>0$. To this end, we consider the following barrier functional structure
\begin{equation} \label{bbbasjka}
B(t,u(t,x)) = \int_0^1 \left[\begin{smallmatrix} u(t,\theta) \\ \partial_\theta u(t,\theta)  \end{smallmatrix}\right]^\prime   M(\theta) \left[\begin{smallmatrix} u(t,\theta) \\ \partial_\theta u(t,\theta)  \end{smallmatrix}\right] ~\mathrm{d}\theta,
\end{equation}
where $M(\theta) \in \mathbb{R}^{2\times 2}$. Applying Corollary~\ref{cor:fwd} and performing a line search for $\lambda$, the maximum parameter $\lambda$, for which the solutions are $ \mathcal{Y}_u$-safe, is found to be $\lambda = 1.196\pi^2$, for which the barrier functional~\eqref{bbbasjka} was constructed with a degree-16 $M(\theta)$ as given in Appendix B. This is consistent with the numerical experiments  shown in Figure~\ref{fig2}, where the $\mathcal{H}^1_\Omega$-norm of the solution to PDE~\eqref{eq:examplemain} with $\lambda=1.2\pi^2$ was computed  for four different initial conditions $u_0(x) \in  \mathcal{U}_0$ as in~\eqref{eq:u00}.

\section{Conclusion and Future Work}

We developed a method based on barrier certificates for verifying whether  the solutions of a PDE are safe with respect to an unsafe set. Numerical examples illustrated the computation of barrier functional certificates by SDPs for problems with polynomial data and equations in one-dimensional spatial domain.

Prospective research can consider bounding functionals of the states of nonlinear stochastic differential equations (SDEs)~\cite{AHP16}, where a method for safety verification of backward-in-time PDEs is developed and used to bound state functionals of SDEs thanks to the Feynman-Kac PDE. This method also has direct applications to optimal  control of stochastic systems, wherein the Hamilton-Jacobi-Bellman equation can be used.

\section*{References}

\bibliography{refsChap4.bib}


\appendix

\section{Well-posedness of PDE Systems} \label{app:wellposed}

We briefly review some aspects related to the well-posedness of PDEs. In the case where $\mathscr{F}$ is a linear operator, the well-posedness problem of (1) is tied to whether $\mathscr{F}$ generates a strongly continuous semigroup denoted $C^0$- Semigroup~\cite[Chapter 2.1]{CZ95}. In this respect, the Hille-Yosida theorem~\cite[Theorem 3.4.1]{Sta05}, \cite[Theorem 2.1.12]{CZ95} provides necessary and sufficient conditions for such generators. In addition, given an operator, the Lumer-Phillips theorem~\cite[ Theorem 3.4.5]{Sta05},~\cite{lumer1961},~\cite[ Theorem 3.8.6]{lumer1961} presents conditions for the generator of a strongly continuous semigroup that are easier to verify based on checking whether the operator is dissipative.

If $\mathscr{F}$ is a nonlinear dissipative operator satisfying $$Dom(\mathscr{F}) \subset Ran(\mathscr{I} - \lambda \mathscr{F}),~\forall \lambda > 0,$$
with $\mathscr{I}$ representing the identity operator, then $\mathscr{F}$ generates a (nonlinear) semigroup of contractions~\cite[Corollary 2.10]{Mi92}. In addition, uniqueness and existence of the solutions to (1) follows from \cite[Theorem 4.10 and Theorem 5.1]{Mi92}.

\section{Numerical Results}

Neglecting the terms with coefficients smaller than $10^{-4}$, the constructed certificate for Example 1 is given by 
\begin{align}
10^4b(t,\theta) =& -7.916\theta^6 + 105.7\theta^5t + 195.0\theta^5 - 315.15\theta^4t^2\nonumber \\
&+ 175.7\theta^4t - 348.2\theta^4 - 35.99\theta^3t^3 - 26.33\theta^3t\nonumber \\
&- 72.06\theta^3 + 42.64\theta^2t^3 - 66.52\theta^2t^2 + 203.8\theta^2t\nonumber \\
& - 228.9\theta^2 - 2.782\theta t^5 - 4.065\theta t^4 - 228.9\theta^2\nonumber \\
&- 2.782\theta t^5 - 4.065\theta t^4 - 1.184\theta t^2 + 2.485\theta t\nonumber\\
& - 15.97\theta - 631.9t^6 + 62.17t^5 - 162.0t^4\nonumber \\
&+ 230.8t^3 - 59.17t^2 + 717.7t - 705.7. \nonumber
\end{align}

Neglecting the terms with coefficients smaller than $10^{-4}$, the constructed certificate for Example 2 is given by

\begin{equation*}
M(\theta) = \left[\begin{matrix} M_{11}(\theta) & M_{12}(x) \\ M_{12}(\theta)  & M_{22}(\theta)   \end{matrix}\right],
\end{equation*}
\begin{align*}
10^{4}M_{11}(\theta) =& - 12.96\theta^{16} + 27.92\theta^{15} - 55.38\theta^{14} - 160.6\theta^{13}
 - 222.4\theta^{12} + 180.8\theta^{11} \\
 &+ 199.1\theta^{10} + 332.9\theta^9
  - 343.5\theta^8 - 454.9\theta^7 - 390.1\theta^6 + 329.9\theta^5\\
 &+ 666.7\theta^4 - 83.37\theta^3 - 663.4\theta^2 + 418.7\theta - 74.97, \nonumber
\end{align*}
\vspace{-3mm}
\begin{align*}
10^{4}M_{12}(\theta) =&  1.39\theta^{16} - 26.03\theta^{15} + 10.76\theta^{14} 
+ 22.53\theta^{13} - 14.63\theta^{12} - 22.81\theta^{11} \\&+ 52.28\theta^{10}
 - 67.56\theta^9 - 69.45\theta^8 - 87.54\theta^7 + 79.37\theta^6 
 + 262.8\theta^5 \\&- 32.63\theta^4 - 447.1\theta^3 + 417.7\theta^2 
 - 157.6\theta + 23.88, \nonumber
 \end{align*}
 \begin{align*}
10^{4}M_{22}(\theta) =&  - 1.607\theta^{16}  - 26.85\theta^{14} + 47.17\theta^{13} 
+ 38.69\theta^{12} - 77.1\theta^{11} - 34.36\theta^{10} \\&+ 66.47\theta^9 
+ 13.36\theta^8 - 34.57\theta^7 - 1.477\theta^6 + 17.13\theta^5 \\
&- 9.405\theta^4 + 2.768\theta^3. \nonumber
 \end{align*}

\end{document}